\documentclass[a4paper,12pt,oneside]{article}
\usepackage{amsmath,amsfonts,amssymb,amsmath,latexsym,amsthm,enumerate}

\newtheorem{thm}{Theorem}
\newtheorem{cor}[thm]{Corollary}
\theoremstyle{definition}

\theoremstyle{plain}

\begin{document}
\title {Poly-Bernoulli polynomials arising from umbral calculus}
\author{by \\Dae San Kim, Taekyun Kim and Sang-Hun Lee}\date{}\maketitle

\begin{abstract}
\noindent In this paper, we give some recurrence formula and new and interesting identities for the poly-Bernoulli numbers and polynomials which are derived from umbral calculus.
\end{abstract}

\section{Introduction}

The classical polylogarithmic function $Li_{s}(x)$ are
\begin{equation}\label{eq:1}
Li_{s}(x)=\sum_{k=1}^{\infty}\frac{x^{k}}{k^{s}},\,\,\,\,s\in\mathbf{Z},\,\,\,\,(\text{see}\,\, \lbrack 3,5\rbrack).
\end{equation}
In $\lbrack 5\rbrack$, poly-Bernoulli polynomials are defined by the generating function to be
\begin{equation}\label{eq:2}
\frac{Li_{k}\left(1-e^{-t}\right)}{1-e^{-t}}e^{xt}=e^{B^{(k)}(x)t}=\sum_{n=0}^{\infty}B_{n}^{(k)}(x)\frac{t^{n}}{n!},\,\,\,(\text{see}\,\, \lbrack 3,5\rbrack),
\end{equation}
with the usual convention about replacing $\left(B^{(k)}(x)\right)^{n}$ by $B_{n}^{(k)}(x)$.\\
As is well known, the Bernoulli polynomials of order $r$ are defined by the generating function to be
\begin{equation}\label{eq:3}
\left(\frac{t}{e^{t}-1}\right)^{r}e^{xt}=\sum_{n=0}^{\infty}\mathbb{B}_{n}^{(r)}(x)\frac{t^{n}}{n!},\,\,\,(\text{see}\,\, \lbrack 7,9\rbrack).
\end{equation}
In the special case, $r=1$, $\mathbb{B}_{n}^{(r)}(x)=B_{n}(x)$ is called the $n$-th ordinary Bernoulli polynomial. Here we denote higher-order Bernoulli polynomials as $\mathbb{B}_{n}^{(r)}$ to avoid conflict of notations.\\
If $x=0$, then $B_{n}^{(k)}(0)=B_{n}^{(k)}$ is called the $n$-th poly-Bernoulli number. From (\ref{eq:2}), we note that
\begin{equation}\label{eq:4}
B_{n}^{(k)}(x)=\sum_{l=0}^{n}\binom{n}{l}B_{n-l}^{(k)}x^{l}=\sum_{l=0}^{n}\binom{n}{l}B_{l}^{(k)}x^{n-l}.
\end{equation}
Let $\mathcal{F}$ be the set of all formal power series in the variable $t$ over $\mathbf{C}$ as follows:
\begin{equation}\label{eq:5}
\mathcal{F}=\left\{f(t)=\sum_{k=0}^{\infty}a_{k}\frac{t^{k}}{k!}\Bigg\vert a_{k}\in\mathbf{C}\right\},
\end{equation}
and let $\mathbb{P}=\mathbf{C}\lbrack x\rbrack$ and $\mathbb{P}^{*}$ denote the vector space of all linear functionals on $\mathbb{P}$. $\left\langle L\vert p(x)\right\rangle$ denotes the acition of linear functional $L$ on the polynomial $p(x)$, and it is well known that the vector space oprations on $\mathbb{P}^{*}$ are defined by $\left\langle L+M\vert p(x)\right\rangle=\left\langle L\vert p(x)\right\rangle+\left\langle M\vert p(x)\right\rangle$, $\left\langle cL\vert p(x)\right\rangle=c\left\langle L\vert p(x)\right\rangle$, where $c$ is a complex constant (see $\lbrack 6,9\rbrack$).\\
For $f(t)\in\mathcal{F}$, let $\left\langle f(t)\big\vert x^{n}\right\rangle=a_{n}$. Then, by (\ref{eq:5}), we easily get
\begin{equation}\label{eq:6}
\left\langle t^{k}\big\vert x^{n}\right\rangle=n!\delta_{n,k},\,\,(n,k\geq 0),\,\,\,(\text{see}\,\, \lbrack 1,4,6,9,10\rbrack),
\end{equation}
where $\delta_{n,k}$ is the Kronecker's symbol.\\
Let us assume that $f_{L}(t)=\sum_{k=0}^{\infty}\left\langle L\big\vert x^{k}\right\rangle\frac{t^{k}}{k!}$. Then, by (\ref{eq:6}), we see that $\left\langle f_{L}(t)\big\vert x^{n}\right\rangle=\left\langle L\big\vert x^{n}\right\rangle$. That is, $f_{L}(t)=L$. Additionally, the map $L\longmapsto f_{L}(t)$ is a vector space isomorphism from $\mathbb{P}^{*}$ onto $\mathcal{F}$. Henceforth, $\mathcal{F}$ denotes both the algebra of the formal power series in $t$ and the vector space of all linear functionals on $\mathbb{P}$, and so an element $f(t)$ of $\mathcal{F}$ will be thought as a formal power series and a linear functional. $\mathcal{F}$ is called the umbral algebra. The umbral calculus is the study of umbral algebra.
The order $O\left(f(t) \right)$ of the power series $f(t)\neq0$ is the smallest integer  for which $a_{k}$ does not vanish. If $O(f(t))=0$, then $f(t)$ is called an invertible series. If $O(f(t))=1$, then $f(t)$ is called a delta series. For $f(t),g(t)\in\mathcal{F}$, we have
\begin{equation}\label{eq:7}
\left\langle f(t)g(t)\vert p(x)\right\rangle=\left\langle f(t)\vert g(t)p(x)\right\rangle=\left\langle g(t)\vert f(t)p(x)\right\rangle.
\end{equation}
Let $f(t)\in\mathcal{F}$ and $p(x)\in\mathbb{P}$. Then we have
\begin{equation}\label{eq:8}
f(t)=\sum_{k=0}^{\infty}\left\langle f(t)\big\vert x^{k}\right\rangle\frac{t^{k}}{k!},\quad p(x)=\sum_{k=0}^{\infty}\left\langle t^{k}\big\vert p(x)\right\rangle\frac{x^{k}}{k!},\,\,\,\,(\text{see}\,\,\lbrack 6,9\rbrack).
\end{equation}
From (\ref{eq:8}), we can easily derive
\begin{equation}\label{eq:9}
p^{(k)}(x)=\frac{d^{k}p(x)}{dx^{k}}=\sum_{l=k}^{\infty}\frac{\left\langle t^{l}\big\vert p(x)\right\rangle}{(l-k)!}x^{l-k}.
\end{equation}
Thus, by (\ref{eq:8}) and (\ref{eq:9}), we get
\begin{equation}\label{eq:10}
p^{(k)}(0)=\left\langle t^{k}\big\vert p(x)\right\rangle=\left\langle 1\big\vert p^{(k)}(x)\right\rangle.
\end{equation}
Hence, from (\ref{eq:10}), we have
\begin{equation}\label{eq:11}
t^{k}p(x)=p^{(k)}(x)=\frac{d^{k}p(x)}{dx^{k}},\,\,\,\,(\text{see}\,\,\lbrack 1,6,9\rbrack).
\end{equation}
It is easy to show that
\begin{equation}\label{eq:12}
e^{yt}p(x)=p(x+y),\quad\left\langle e^{yt}\big\vert p(x)\right\rangle=p(y).
\end{equation}
Let $O(f(t))=1$ and $O(g(t))=0$. Then there exists a unique sequence $s_{n}(x)$ of polynomials such that $\left\langle g(t)f(t)^{k}\big\vert s_{n}(x)\right\rangle=n!\delta_{n,k}$, for $n, k\geq 0$. The sequence $s_{n}(x)$ is called a Sheffer sequence for $\left(g(t),f(t)\right)$ which is denoted by $s_{n}(x)\sim\left(g(t),f(t)\right)$. The Sheffer sequence $s_{n}(x)$ for $\left(g(t),t\right)$ is called the Appell sequence for $g(t)$. For $p(x)\in\mathbb{P}$, $f(t)\in\mathcal{F}$, we have
\begin{equation}\label{eq:13}
\left\langle f(t)\vert xp(x)\right\rangle=\left\langle \partial_{t}f(t)\vert p(x)\right\rangle=\left\langle f'(t)\vert p(x)\right\rangle,\,\,\,\,(\text{see}\,\,\lbrack 9,10\rbrack).
\end{equation}
Let $s_{n}(x)\sim\left(g(t),f(t)\right)$. Then the following equations are known:
\begin{equation}\label{eq:14}
h(t)=\sum_{k=0}^{\infty}\frac{\left\langle h(t)\big\vert s_{k}(x)\right\rangle}{k!}g(t)f(t)^{k},\quad p(x)=\sum_{k=0}^{\infty}\frac{\left\langle g(t)f(t)^{k}\big\vert p(x)\right\rangle}{k!}s_{k}(x),
\end{equation}
where $h(t)\in\mathcal{F}$, $p(x)\in\mathbb{P}$,
\begin{equation}\label{eq:15}
\frac{1}{g\left(\bar{f}(t)\right)}e^{y\bar{f}(t)}=\sum_{k=0}^{\infty}s_{k}(y)\frac{t^{k}}{k!},\,\,\,\,\text{for all}\,\, y\in\mathbf{C},
\end{equation}
where $\bar{f}(t)$ is the compositional inverse for $f(t)$ with $\bar{f}\left(f(t)\right)=t$, and
\begin{equation}\label{eq:16}
f(t)s_{n}(x)=ns_{n-1}(x).
\end{equation}
As is well known, the Stirling numbers of the second kind are also defined by the generating function to be
\begin{equation}\label{eq:17}
\left(e^{t}-1\right)^{m}=m!\sum_{l=m}^{\infty}S_{2}(l,m)\frac{t^{l}}{l!}=\sum_{l=0}^{\infty}\frac{m!}{(l+m)!}S_{2}(l+m,m)t^{l+m}.
\end{equation}
Let $s_{n}(x)\sim\left(g(t),t\right)$. Then the Appell identity is given by
\begin{equation}\label{eq:18}
s_{n}(x+y)=\sum_{k=0}^{n}\binom{n}{k}s_{k}(y)x^{n-k}=\sum_{k=0}^{n}\binom{n}{k}s_{n-k}(y)x^{k},\,\,\,\,(\text{see}\,\,\lbrack 6,9\rbrack),
\end{equation}
and
\begin{equation}\label{eq:19}
s_{n+1}(x)=\left(x-\frac{g'(t)}{g(t)}\right)s_{n}(x),\,\,\,\,(\text{see}\,\,\lbrack 6.9\rbrack).
\end{equation}
For $s_{n}(x)\sim\left(g(t),f(t)\right)$, $r_{n}(x)\sim\left(h(t),l(t)\right)$, we have
\begin{equation}\label{eq:20}
s_{n}(x)=\sum_{m=0}^{n}r_{m}(x)c_{n,m},
\end{equation}
where
\begin{equation}\label{eq:21}
c_{n,m}=\frac{1}{m!}\left\langle \frac{h\left(\bar{f}(t)\right)}{g\left(\bar{f}(t)\right)}l\left(\bar{f}(t)\right)^{m}\Bigg\vert x^{n}\right\rangle,\,\,\,\,(\text{see}\,\,\lbrack 9\rbrack).
\end{equation}
The equations (\ref{eq:20}) and (\ref{eq:21}) are important in deriving our main results of this paper.\\
In this paper, we give some recurrence formula and new and interesting identities for the poly-Bernoulli numbers and polynomials which are derive from umbral calculus.


\section{Poly-Bernoulli numbers and polynomials}

Let $g_{k}(t)=\frac{1-e^{-t}}{Li_{k}\left(1-e^{-t}\right)}$. Then, by (\ref{eq:2}) and (\ref{eq:15}), we get
\begin{equation}\label{eq:22}
B_{n}^{(k)}(x)\sim\left(g_{k}(t),t\right).
\end{equation}
That is, poly-Bernoulli polynomial $B_{n}^{(k)}(x)$ is an Appell sequence.\\
By (\ref{eq:1}), we easily get
\begin{equation}\label{eq:23}
\frac{d}{dx}Li_{k}(x)=\frac{1}{x}Li_{k-1}(x),\quad tB_{n}^{(k)}(x)=\frac{d}{dx}B_{n}^{(k)}(x)=nB_{n-1}^{(k)}(x).
\end{equation}
From (\ref{eq:2}) and (\ref{eq:15}), we have
\begin{equation}\label{eq:24}
B_{n}^{(k)}(x)=\frac{1}{g_{k}(t)}x^{n}=\frac{Li_{k}\left(1-e^{-t}\right)}{1-e^{-t}}x^{n}.
\end{equation}
Let $k\in\mathbf{Z}$ and $n\geq 0$. Then we have
\begin{align}\label{eq:25}
B_{n}^{(k)}(x)&=\frac{Li_{k}\left(1-e^{-t}\right)}{1-e^{-t}}x^{n}=\sum_{m=1}^{\infty}\frac{\left(1-e^{-t}\right)^{m-1}}{m^{k}}x^{n}\\
&=\sum_{m=0}^{\infty}\frac{1}{(m+1)^{k}}\left(1-e^{-t}\right)^{m}x^{n}=\sum_{m=0}^{\infty}\frac{1}{(m+1)^{k}}\sum_{j=0}^{m}(-1)^{j}\binom{m}{j}e^{-jt}x^{n}\nonumber\\
&=\sum_{m=0}^{n}\frac{1}{(m+1)^{k}}\sum_{j=0}^{m}(-1)^{j}\binom{m}{j}(x-j)^{n}.\nonumber
\end{align}
By (\ref{eq:17}) and (\ref{eq:25}), we get
\begin{align}\label{eq:26}
B_{n}^{(k)}(x)&=\sum_{m=0}^{n}\frac{1}{(m+1)^{k}}\sum_{a=0}^{\infty}(-1)^{a}\frac{m!}{(a+m)!}S_{2}(a+m,m)t^{a+m}x^{n}\\
&=\sum_{m=0}^{n}\frac{1}{(m+1)^{k}}\sum_{a=0}^{n-m}(-1)^{a}\frac{m!}{(a+m)!}S_{2}(a+m,m)(n)_{a+m}x^{n-a-m}\nonumber\\
&=\sum_{l=0}^{n}\left\{\sum_{m=0}^{n-l}\frac{(-1)^{n-m-l}}{(m+1)^{k}}\binom{n}{l}m!S_{2}(n-l,m)\right\}x^{l},\nonumber
\end{align}
where $(a)_{n}=a(a-1)(a-2)\cdots (a-n+1)$.\\
Now, we use the well-known transfer formula for Appell sequences (see equation (\ref{eq:19})).\\
By (\ref{eq:19}) and (\ref{eq:22}), we get
\begin{equation}\label{eq:27}
B_{n+1}^{(k)}(x)=\left(x-\frac{g_{k}'(t)}{g_{k}(t)}\right)B_{n}^{(k)}(x),
\end{equation}
where
\begin{align}\label{eq:28}
\frac{g_{k}'(t)}{g_{k}(t)}&=\left(\log{g_{k}(t)}\right)'=\left(\log{\left(1-e^{-t}\right)}-\log{Li_{k}}\left(1-e^{-t}\right)\right)'\\
&=\frac{e^{-t}}{1-e^{-t}}\left\{1-\frac{{Li}_{k-1}\left(1-e^{-t}\right)}{{Li}_{k}\left(1-e^{-t}\right)}\right\}\nonumber\\
&=\frac{1}{e^{t}-1}\left(\frac{Li_{k}\left(1-e^{-t}\right)-Li_{k-1}\left(1-e^{-t}\right)}{Li_{k}\left(1-e^{-t}\right)}\right).\nonumber
\end{align}
From (\ref{eq:27}) and (\ref{eq:28}), we have
\begin{align}\label{eq:29}
B_{n+1}^{(k)}(x)&=xB_{n}^{(k)}(x)-\frac{g_{k}'(t)}{g_{k}(t)}B_{n}^{(k)}(x)\\
&=xB_{n}^{(k)}(x)-\left(\frac{t}{e^{t}-1}\right)\left(\frac{Li_{k}\left(1-e^{-t}\right)-Li_{k-1}\left(1-e^{-t}\right)}{t\left(1-e^{-t}\right)}\right)x^{n}.\nonumber
\end{align}
Here, we note that
\begin{align}\label{eq:30}
\frac{Li_{k}\left(1-e^{-t}\right)-Li_{k-1}\left(1-e^{-t}\right)}{1-e^{-t}}&=\sum_{m=2}^{\infty}\left(\frac{1}{m^{k}}-\frac{1}{m^{k-1}}\right)\left(1-e^{-t}\right)^{m-1}\\
&=\left(\frac{1}{2^{k}}-\frac{1}{2^{k-1}}\right)t+\cdots\nonumber
\end{align}
is a delta series.\\
For any delta series $f(t)$, we observe that
\begin{equation}\label{eq:31}
\frac{f(t)}{t}x^{n}=f(t)\frac{x^{n+1}}{n+1}.
\end{equation}
By (\ref{eq:29}), (\ref{eq:30}) and (\ref{eq:31}), we get
\begin{align}\label{eq:32}
B_{n+1}^{(k)}(x)&=xB_{n}^{(k)}(x)-\left(\frac{t}{e^{t}-1}\right)\left(\frac{1}{n+1}\frac{Li_{k}\left(1-e^{-t}\right)-Li_{k-1}\left(1-e^{-t}\right)}{1-e^{-t}}x^{n+1}\right)\nonumber\\
&=xB_{n}^{(k)}(x)-\frac{1}{n+1}\sum_{l=0}^{\infty}\frac{B_{l}}{l!}t^{l}\left\{B_{n+1}^{(k)}(x)-B_{n+1}^{(k-1)}(x)\right\}\nonumber\\
&=xB_{n}^{(k)}(x)-\frac{1}{n+1}\sum_{l=0}^{n+1}\binom{n+1}{l}B_{l}\left\{B_{n+1-l}^{(k)}(x)-B_{n+1-l}^{(k-1)}(x)\right\}.
\end{align}
Therefore, by (\ref{eq:32}), we obtain the following theorem.

\begin{thm}\label{eq:thm1}
For $k\in\mathbf{Z}$, $n\geq 0$, we have
\begin{equation*}
B_{n+1}^{(k)}(x)=xB_{n}^{(k)}(x)-\frac{1}{n+1}\sum_{l=0}^{n+1}\binom{n+1}{l}B_{l}\left\{B_{n+1-l}^{(k)}(x)-B_{n+1-l}^{(k-1)}(x)\right\},
\end{equation*}
where $B_{n}$ is the $n$-th ordinary Bernoulli number.
\end{thm}

\noindent It is easy to show that
\begin{align}\label{eq:33}
txB_{n}^{(k)}(x)&=t\sum_{l=0}^{n}\binom{n}{l}B_{n-l}^{(k)}x^{l+1}=\sum_{l=0}^{n}\binom{n}{l}B_{n-l}^{(k)}(l+1)x^{l}\\
&=nx\sum_{l=0}^{n-1}\binom{n-1}{l}B_{n-1-l}^{(k)}x^{l}+\sum_{l=0}^{n}\binom{n}{l}B_{n-l}^{(k)}x^{n}\nonumber\\
&=nxB_{n-1}^{(k)}(x)+B_{n}^{(k)}(x).\nonumber
\end{align}
Applying to $t$ on both sides of Theorem \ref{eq:thm1}, by (\ref{eq:33}), we get
\begin{equation}\label{eq:34}
(n+1)B_{n}^{(k)}(x)=nxB_{n-1}^{(k)}(x)+B_{n}^{(k)}(x)-\sum_{l=0}^{n}\binom{n}{l}B_{n-l}\left\{B_{l}^{(k)}(x)-B_{l}^{(k-1)}(x)\right\}.
\end{equation}
Therefore, by (\ref{eq:34}), we obtain the following corollary.

\begin{cor}\label{eq:cor2}
For $k\in\mathbf{Z}$ and $n\geq 1$, we have
\begin{align*}
&(n+1)B_{n}^{(k)}(x)-n\left(x+\frac{1}{2}\right)B_{n-1}^{(k)}(x)+\sum_{l=0}^{n-2}\binom{n}{l}B_{n-l}B_{l}^{(k)}(x)\\
&=\sum_{l=0}^{n}\binom{n}{l}B_{n-l}B_{l}^{(k-1)}(x).
\end{align*}
\end{cor}

\noindent From (\ref{eq:2}) and (\ref{eq:6}), we note that
\begin{align}\label{eq:35}
B_{n}^{(k)}(y)&=\left\langle\frac{Li_{k}\left(1-e^{-t}\right)}{1-e^{-t}}e^{yt}\Bigg\vert x^{n}\right\rangle=\left\langle\frac{Li_{k}\left(1-e^{-t}\right)}{1-e^{-t}}e^{yt}\Bigg\vert xx^{n-1}\right\rangle\\
&=\left\langle\partial_{t}\left(\frac{Li_{k}\left(1-e^{-t}\right)}{1-e^{-t}}e^{yt}\right)\Bigg\vert x^{n-1}\right\rangle\nonumber\\
&=\left\langle\partial_{t}\left(\frac{Li_{k}\left(1-e^{-t}\right)}{1-e^{-t}}\right)e^{yt}\Bigg\vert x^{n-1}\right\rangle +y\left\langle\frac{Li_{k}\left(1-e^{-t}\right)}{1-e^{-t}}e^{yt}\Bigg\vert x^{n-1}\right\rangle\nonumber\\
&=\left\langle\frac{Li_{k-1}\left(1-e^{-t}\right)-Li_{k}\left(1-e^{-t}\right)}{\left(1-e^{-t}\right)^{2}}e^{(y-1)t}\Bigg\vert x^{n-1}\right\rangle +yB_{n-1}^{(k)}(y).\nonumber
\end{align}
Now, we observe that
\begin{align}\label{eq:36}
\frac{Li_{k-1}\left(1-e^{-t}\right)-Li_{k}\left(1-e^{-t}\right)}{\left(1-e^{-t}\right)^{2}}&=\frac{1}{\left(1-e^{-t}\right)^{2}}\sum_{m=1}^{\infty}\left\{\frac{\left(1-e^{-t}\right)^{m}}{m^{k-1}}-\frac{\left(1-e^{-t}\right)^{m}}{m^{k}}\right\}\\
&=\sum_{m=2}^{\infty}\left\{\frac{1}{m^{k-1}}-\frac{1}{m^{k}}\right\}\left(1-e^{-t}\right)^{m-2}\nonumber\\
&=\sum_{m=0}^{\infty}\left\{\frac{1}{(m+2)^{k-1}}-\frac{1}{(m+2)^{k}}\right\}\left(1-e^{-t}\right)^{m}.\nonumber
\end{align}
Thus, by (\ref{eq:36}), we get
\begin{align}\label{eq:37}
&\left\langle\frac{Li_{k-1}\left(1-e^{-t}\right)-Li_{k}\left(1-e^{-t}\right)}{\left(1-e^{-t}\right)^{2}}e^{(y-1)t}\Bigg\vert x^{n-1}\right\rangle\\
&=\sum_{m=0}^{\infty}\left\{\frac{1}{(m+2)^{k-1}}-\frac{1}{(m+2)^{k}}\right\}\left\langle\left(1-e^{-t}\right)^{m}e^{(y-1)t}\big\vert x^{n-1}\right\rangle\nonumber\\
&=\sum_{m=0}^{n-1}\left(\frac{1}{(m+2)^{k-1}}-\frac{1}{(m+2)^{k}}\right)\left\langle\left(1-e^{-t}\right)^{m}\big\vert\left(x+y-1\right)^{n-1}\right\rangle\nonumber\\
&=\sum_{m=0}^{n-1}\left(\frac{1}{(m+2)^{k-1}}-\frac{1}{(m+2)^{k}}\right)\sum_{a=0}^{n-1}\binom{n-1}{a}(y-1)^{n-1-a}\left\langle\left(1-e^{-t}\right)^{m}\big\vert x^{a}\right\rangle.\nonumber
\end{align}
From (\ref{eq:6}) and (\ref{eq:7}), we have
\begin{equation}\label{eq:38}
\left\langle\left(1-e^{-t}\right)^{m}\big\vert x^{a}\right\rangle=(-1)^{a+m}m!S_{2}(a,m).
\end{equation}
From (\ref{eq:37}) and (\ref{eq:38}), wehave
\begin{align}\label{eq:39}
&\left\langle\frac{Li_{k-1}\left(1-e^{-t}\right)-Li_{k-1}\left(1-e^{-t}\right)}{\left(1-e^{-t}\right)^{2}}e^{(y-1)t}\Bigg\vert x^{n-1}\right\rangle\\
&=\sum_{m=0}^{n-1}\sum_{a=0}^{n-1}(-1)^{a+m}\binom{n-1}{a}m!\left(\frac{1}{(m+2)^{k-1}}-\frac{1}{(m+2)^{k}}\right)\nonumber\\
&\quad\times S_{2}(a,m)(y-1)^{n-1-a}\nonumber\\
&=\sum_{m=0}^{n-1}\sum_{l=0}^{n-1}(-1)^{n-1-l+m}\binom{n-1}{l}m!\left(\frac{1}{(m+2)^{k-1}}-\frac{1}{(m+2)^{k}}\right)\nonumber\\
&\quad\times S_{2}(n-1-l,m)(y-1)^{l}.\nonumber
\end{align}
Therefore, by (\ref{eq:35}) and (\ref{eq:39}), we obtain the following theorem.

\begin{thm}\label{eq:thm3}
For $k\in\mathbf{Z}$ and $n\geq 1$, we have
\begin{align*}
B_{n}^{(k)}(x)&=xB_{n-1}^{(k)}(x)+\sum_{l=0}^{n-1}(-1)^{n-1-l}\binom{n-1}{l}\\
&\quad\times\left\{\sum_{m=0}^{n-1}(-1)^{m}\frac{(m+1)!}{(m+2)^{k}}S_{2}(n-1-l,m)\right\}(x-1)^{l}.
\end{align*}
\end{thm}
\noindent Now, we try to compute $\left\langle Li_{k}\left(1-e^{-t}\right)\big\vert x^{n+1}\right\rangle$ in two ways. On the one hand,
\begin{align}\label{eq:40}
\left\langle Li_{k}\left(1-e^{-t}\right)\big\vert x^{n+1}\right\rangle&=\left\langle \left(1-e^{-t}\right)\frac{Li_{k}\left(1-e^{-t}\right)}{1-e^{-t}}\Bigg\vert x^{n+1}\right\rangle\\
&=\left\langle\frac{Li_{k}\left(1-e^{-t}\right)}{1-e^{-t}}\Bigg\vert \left(1-e^{-t}\right)x^{n+1}\right\rangle\nonumber\\
&=\left\langle\frac{Li_{k}\left(1-e^{-t}\right)}{1-e^{-t}}\Bigg\vert x^{n+1}-(x-1)^{n+1}\right\rangle\nonumber\\
&=\sum_{m=0}^{n}\binom{n+1}{m}(-1)^{n-m}\left\langle 1\Bigg\vert \frac{Li_{k}\left(1-e^{-t}\right)}{1-e^{-t}}x^{m}\right\rangle\nonumber\\
&=\sum_{m=0}^{n}\binom{n+1}{m}(-1)^{n-m}B_{m}^{(k)}.\nonumber
\end{align}
On the other hand,
\begin{align}\label{eq:41}
&\left\langle Li_{k}\left(1-e^{-t}\right)\big\vert x^{n+1}\right\rangle=\left\langle\int_{0}^{t}\left(Li_{k}\left(1-e^{-s}\right)\right)'ds\Bigg\vert x^{n+1}\right\rangle\\
&=\left\langle\int_{0}^{t}e^{-s}\frac{Li_{k-1}\left(1-e^{-s}\right)}{1-e^{-s}}ds\Bigg\vert x^{n+1}\right\rangle\nonumber\\
&=\left\langle\int_{0}^{t}\left(\sum_{a=0}^{\infty}\frac{(-s)^{a}}{a!}\right)\left(\sum_{m=0}^{\infty}\frac{B_{m}^{(k-1)}}{m!}s^{m}\right)ds\Bigg\vert x^{n+1}\right\rangle\nonumber\\
&=\left\langle\sum_{l=0}^{\infty}\left(\sum_{m=0}^{l}\binom{l}{m}(-1)^{l-m}B_{m}^{(k-1)}\right)\frac{1}{l!}\int_{0}^{t}s^{l}ds\Bigg\vert x^{n+1}\right\rangle\nonumber\\
&=\sum_{l=0}^{n}\sum_{m=0}^{l}\binom{l}{m}(-1)^{l-m}\frac{B_{m}^{(k-1)}}{(l+1)!}(n+1)!\delta_{n+1,l+1}=\sum_{m=0}^{n}\binom{n}{m}(-1)^{n-m}B_{m}^{(k-1)}.\nonumber
\end{align}
By (\ref{eq:40}) and (\ref{eq:41}), we get
\begin{equation}\label{eq:42}
\sum_{m=0}^{n}(-1)^{n-m}\binom{n}{m}B_{m}^{(k-1)}=\sum_{m=0}^{n}\binom{n+1}{m}(-1)^{n-m}B_{m}^{(k)}.
\end{equation}
By (\ref{eq:2}) and (\ref{eq:3}), we see that
\begin{equation}\label{eq:43}
B_{n}^{(k)}(x)\sim\left(\frac{1-e^{-t}}{Li_{k}\left(1-e^{-t}\right)},t\right),\quad\mathbb{B}_{n}^{(r)}(x)\sim\left(\left(\frac{e^{t}-1}{t}\right)^{r},t\right),\,\,\,\, r\geq 0.
\end{equation}
From (\ref{eq:20}), (\ref{eq:21}) and (\ref{eq:43}), we have
\begin{equation}\label{eq:44}
B_{n}^{(k)}(x)=\sum_{m=0}^{n}C_{n,m}\mathbb{B}_{m}^{(r)}(x),
\end{equation}
where
\begin{align}\label{eq:45}
C_{n,m}&=\frac{1}{m!}\left\langle\frac{\left(\frac{e^{t}-1}{t}\right)^{r}}{\frac{1-e^{-t}}{Li_{k}\left(1-e^{-t}\right)}}t^{m}\Bigg\vert x^{n}\right\rangle\\
&=\frac{1}{m!}\left\langle\frac{Li_{k}\left(1-e^{-t}\right)}{1-e^{-t}}\left(\frac{e^{t}-1}{t}\right)^{r}\Bigg\vert t^{m}x^{n}\right\rangle\nonumber\\
&=\binom{n}{m}\left\langle\frac{Li_{k}\left(1-e^{-t}\right)}{1-e^{-t}}\Bigg\vert\left(\frac{e^{t}-1}{t}\right)^{r}x^{n-m}\right\rangle.\nonumber
\end{align}
By (\ref{eq:17}), we easily get
\begin{equation}\label{eq:46}
\left(\frac{e^{t}-1}{t}\right)^{r}=\sum_{l=0}^{\infty}\frac{r!}{(l+r)!}S_{2}(l+r,r)t^{l}.
\end{equation}
Thus, from (\ref{eq:46}), we have
\begin{equation}\label{eq:47}
\left(\frac{e^{t}-1}{t}\right)^{r}x^{n-m}=\sum_{l=0}^{n-m}\frac{r!}{(l+r)!}S_{2}(l+r,r)(n-m)_{l}x^{n-m-l}.
\end{equation}
By (\ref{eq:45}) and (\ref{eq:47}), we get
\begin{align}\label{eq:48}
C_{n,m}&=\binom{n}{m}\sum_{l=0}^{n-m}\frac{r!}{(l+r)!}S_{2}(l+r,r)(n-m)_{l}\left\langle\frac{Li_{k}\left(1-e^{-t}\right)}{1-e^{-t}}\Bigg\vert x^{n-m-l}\right\rangle\\
&=\binom{n}{m}\sum_{l=0}^{n-m}\frac{r!}{(l+r)!}S_{2}(l+r,r)(n-m)_{l}\left\langle t^{0}\Bigg\vert \frac{Li_{k}\left(1-e^{-t}\right)}{1-e^{-t}}x^{n-m-l}\right\rangle\nonumber\\
&=\binom{n}{m}\sum_{l=0}^{n-m}\frac{r!}{(l+r)!}S_{2}(l+r,r)(n-m)_{l}B_{n-m-l}^{(k)}.\nonumber
\end{align}
Therefore, by (\ref{eq:44}) and (\ref{eq:48}), we obtain the following theorem.

\begin{thm}\label{eq:thm4}
For $k\in\mathbf{Z}$ and $r\in\mathbf{Z}_{\geq 0}$, we have
\begin{align*}
B_{n}^{(k)}(x)=\sum_{m=0}^{n}\left\{\binom{n}{m}\sum_{l=0}^{n-m}\frac{r!(n-m)_{l}}{(l+r)!}S_{2}(l+r,r)B_{n-l-m}^{(k)}\right\}\mathbb{B}_{n}^{(r)}(x).
\end{align*}
\end{thm}

\noindent For $r\in\mathbf{Z}_{\geq 0}$, the Euler polynomials of order $r$ are defined by the generating function to be
\begin{equation}\label{eq:49}
\left(\frac{2}{e^{t}+1}\right)^{r}e^{xt}=\sum_{n=0}^{\infty}E_{n}^{(r)}(x)\frac{t^{n}}{n!},\,\,\,\,(\text{see}\,\, \lbrack 2,6,8\rbrack).
\end{equation}
By (\ref{eq:2}) and (\ref{eq:49}), we see that
\begin{equation}\label{eq:50}
B_{n}^{(k)}(x)\sim\left(\frac{1-e^{-t}}{Li_{k}\left(1-e^{-t}\right)},t\right),\quad E_{n}^{(r)}(x)\sim\left(\left(\frac{e^{t}+1}{2}\right)^{r},t\right).
\end{equation}
From (\ref{eq:20}), (\ref{eq:21}) and (\ref{eq:50}), we have
\begin{equation}\label{eq:51}
B_{n}^{(k)}(x)=\sum_{m=0}^{n}C_{n,m}E_{m}^{(r)}(x),
\end{equation}
where
\begin{align}\label{eq:52}
C_{n,m}&=\frac{1}{m!}\left\langle\frac{Li_{k}\left(1-e^{-t}\right)}{1-e^{-t}}\left(\frac{e^{t}+1}{2}\right)^{r}\Bigg\vert t^{m}x^{n}\right\rangle\\
&=\frac{\binom{n}{m}}{2^{r}}\left\langle\frac{Li_{k}\left(1-e^{-t}\right)}{1-e^{-t}}\Bigg\vert \left(e^{t}+1\right)^{r}x^{n-m}\right\rangle\nonumber\\
&=\frac{\binom{n}{m}}{2^{r}}\sum_{j=0}^{r}\binom{r}{j}\left\langle\frac{Li_{k}\left(1-e^{-t}\right)}{1-e^{-t}}\Bigg\vert e^{jt}x^{n-m}\right\rangle\nonumber
\end{align}
\begin{align*}
&=\frac{\binom{n}{m}}{2^{r}}\sum_{j=0}^{r}\binom{r}{j}\left\langle t^{0}\Bigg\vert\frac{Li_{k}\left(1-e^{-t}\right)}{1-e^{-t}}(x+j)^{n-m}\right\rangle\nonumber\\
&=\frac{\binom{n}{m}}{2^{r}}\sum_{j=0}^{r}\binom{r}{j}B_{n-m}^{(k)}(j).
\end{align*}
Therefore, by (\ref{eq:51}) and (\ref{eq:52}), we obtain the following theorem.

\begin{thm}\label{eq:thm5}
For $k\in\mathbf{Z}$ and $r\in\mathbf{Z}_{\geq 0}$, we have
\begin{align*}
B_{n}^{(k)}(x)=\frac{1}{2^{r}}\sum_{m=0}^{n}\left\{\binom{n}{m}\sum_{j=0}^{r}\binom{r}{j}B_{n-m}^{(k)}(j)\right\}E_{n}^{(r)}(x).
\end{align*}
\end{thm}

\noindent Let $\lambda\in\mathbf{C}$ with $\lambda\neq 1$. For $r\in\mathbf{Z}_{\geq 0}$, the Frobenius-Euler polynomials are also defined by the generating function to be
\begin{equation}\label{eq:53}
\left(\frac{1-\lambda}{e^{t}-\lambda}\right)^{r}e^{xt}=\sum_{n=0}^{\infty}H_{n}^{(r)}(x\vert\lambda)\frac{t^{n}}{n!},\,\,\,\,(\text{see}\,\, \lbrack 2,6,7\rbrack).
\end{equation}
From (\ref{eq:2}), (\ref{eq:15}) and (\ref{eq:53}), we note that
\begin{equation}\label{eq:54}
B_{n}^{(k)}(x)\sim\left(\frac{1-e^{-t}}{Li_{k}\left(1-e^{-t}\right)},t\right),\quad H_{n}^{(r)}(x\vert\lambda)\sim\left(\left(\frac{e^{t}-\lambda}{1-\lambda}\right)^{r},t\right).
\end{equation}
By (\ref{eq:20}), (\ref{eq:21}) and (\ref{eq:54}), we get
\begin{equation}\label{eq:55}
B_{n}^{(k)}(x)=\sum_{m=0}^{n}C_{n,m}H_{m}^{(r)}(x\vert\lambda),
\end{equation}
where
\begin{align}\label{eq:56}
C_{n,m}&=\frac{1}{m!}\left\langle\frac{Li_{k}\left(1-e^{-t}\right)}{1-e^{-t}}\left(\frac{e^{t}-\lambda}{1-\lambda}\right)^{r}\Bigg\vert t^{m}x^{n}\right\rangle\\
&=\frac{\binom{n}{m}}{(1-\lambda)^{r}}\left\langle\frac{Li_{k}\left(1-e^{-t}\right)}{1-e^{-t}}\Bigg\vert\left(e^{t}-\lambda\right)^{r}x^{n-m}\right\rangle\nonumber\\
&=\frac{\binom{n}{m}}{(1-\lambda)^{r}}\sum_{j=0}^{r}\binom{r}{j}(-\lambda)^{r-j}\left\langle\frac{Li_{k}\left(1-e^{-t}\right)}{1-e^{-t}}\Bigg\vert e^{jt}x^{n-m}\right\rangle\nonumber\\
&=\frac{\binom{n}{m}}{(1-\lambda)^{r}}\sum_{j=0}^{r}\binom{r}{j}(-\lambda)^{r-j}\left\langle t^{0}\Bigg\vert\frac{Li_{k}\left(1-e^{-t}\right)}{1-e^{-t}}(x+j)^{n-m}\right\rangle\nonumber\\
&=\frac{\binom{n}{m}}{(1-\lambda)^{r}}\sum_{j=0}^{r}\binom{r}{j}(-\lambda)^{r-j}B_{n-m}^{(k)}(j).\nonumber
\end{align}
Therefore, by (\ref{eq:55}) and (\ref{eq:56}), we obtain the following theorem.

\begin{thm}\label{eq:thm6}
For $k\in\mathbf{Z}$ and $r\in\mathbf{Z}_{\geq 0}$, we have
\begin{align*}
B_{n}^{(k)}(x)=\frac{1}{(1-\lambda)^{r}}\sum_{m=0}^{n}\left\{\binom{n}{m}\sum_{j=0}^{r}\binom{r}{j}(-\lambda)^{r-j}B_{n-m}^{(k)}(j)\right\}H_{m}^{(r)}(x\vert\lambda).
\end{align*}
\end{thm}


\noindent
\author{Department of Mathematics, Sogang University, Seoul 121-742, Republic of Korea
\\e-mail: dskim@sogang.ac.kr}\\
\\
\author{Department of Mathematics, Kwangwoon University, Seoul 139-701, Republic of Korea
\\e-mail: tkkim@kw.ac.kr}\\
\author{Division of General Education, Kwangwoon University, Seoul 139-701, Republic of Korea
\\e-mail: leesh58@kw.ac.kr}
\end{document}